\documentclass
[12pt]
{amsart}

\usepackage{amsmath}

\usepackage{amsfonts}
\usepackage{amssymb}

\newtheorem{theorem}{Theorem}[section]
\newtheorem{lemma}[theorem]{Lemma}

\theoremstyle{definition}

\newtheorem{example}[theorem]{Example}

\theoremstyle{remark}
\newtheorem{remark}[theorem]{Remark}

\numberwithin{equation}{section}

\DeclareMathOperator{\R}{{\mathbb{R}}}
\DeclareMathOperator{\T}{{\mathbb{T}}}
\DeclareMathOperator{\Z}{{\mathbb{Z}}}
\DeclareMathOperator{\Int}{{Int}}
\DeclareMathOperator{\Schur}{{Sch}}
\DeclareMathOperator{\Riesz}{{Rsz}}
\DeclareMathOperator{\Alt}{{Alt}}
\DeclareMathOperator{\BohrR}{{\rm b}\!\R}
\DeclareMathOperator{\Rd}{\R_{\rm d}}
\DeclareMathOperator{\Rdtwo}{\R_{\rm d}^2}
\DeclareMathOperator{\Gd}{G_{\rm d}}

\usepackage{caption}
\usepackage{subcaption}
\usepackage{graphicx}

\usepackage{hyperref}

\begin{document}


\title[
Discrete Fourier restriction theorems.
]
{Discrete Fourier restriction theorems\\
in two dimensions.}
\author{John J.F. Fournier}
\address{Department of Mathematics\\
University of British Columbia\\
1984 Mathematics Road\\
Vancouver BC\\
Canada V6T 1Z2}
\email{fournier@math.ubc.ca}
\thanks{
Announced
at the 6th Conference on Function Spaces
in May 2010}
\thanks{Research supported by NSERC grant 4822.}


\begin{abstract}
Consider the group~${\mathbb{R}}^2$ with the discrete topology, and denote its
Fourier algebra by~$A({{\mathbb{R}}_{\rm d}^2})$.
We reformulate a theorem of V.A. Yudin as a statement about restrictions of
functions in~$A({{\mathbb{R}}_{\rm d}^2})$ to
the boundary of a strictly convex domain when those functions vanish outside that boundary.
We give
visual proofs of that statement and
a
complementary one.

\end{abstract}

\maketitle

\markleft{
John J.F. Fournier
}
\pagestyle{myheadings}

\section{Introduction}\label{sec:intro}

Yudin's theorem~\cite{Yudin} is about the Fourier coefficients,~$\hat f(\vec n)$ say,
of an integrable function~$f$
on the product~$\T\times\T$
of two copies of the unit circle group~$\T$.
Those coefficients
are defined on the product~$\Z\times\Z$ of two copies of the integer group~$\Z$.
He used a dual method to estimate the~$\ell^2$ norm of
their
restriction
to the integer lattice points in
the boundary of a strictly convex domain in~$\R^2$
when~$\hat f$ vanishes
outside
that boundary.
We give direct
proofs of that estimate and of the corresponding estimate when~$\hat f$ vanishes
inside
the boundary.

As usual,
\[
\hat f(n_1, n_2) = \left(\frac{1}{2\pi}\right)^2
\int_{-\pi}^\pi
\int_{-\pi}^\pi
f(t_1, t_2)e^{-n_1t_1}e^{-n_2t_2}\,dt_1\,dt_2.
\]
Use the
same
measure~$(1/2\pi)^2dt_1\,dt_2$
in computing~$L^p$ norms.
Given a subset~$D$ of~$\R^2$, denote
its interior by~$\Int(D)$, its complement by~$D^c$ and its boundary by~$\Gamma$.

Our main goal in this paper is to give visual proofs of both parts of an extension of the following statement.

\begin{theorem}
\label{th:Yudin}
There is a constant~$C$ so that
if~$D$ is
a strictly convex set in~$\R^2$ with boundary~$\Gamma$,
and if~$f \in L^1(\T^2)$, then
the estimate
\begin{equation}
\label{eq:YudinIntegers}
\left[\sum_{\vec n \in \Gamma\cap\Z^2} |\hat f(\vec n)|^2\right]^{1/2}
\le C\|f\|_1
\end{equation}
follows from either of the following conditions:
\begin{enumerate}
\item{}
\label{it:interior}
$\hat f$ vanishes
on~$\Int(D)\cap\Z^2$.
\item{}
\label{it:exterior}
$\hat f$ vanishes
on~$\Int(D^c)\cap\Z^2$.
\end{enumerate}
\end{theorem}

Call these the \emph{interior}
and \emph{exterior} cases.
As in~\cite[p. 861]{Yudin},
no uniform estimate of the form~\eqref{eq:YudinIntegers} is possible
in either case
for a
family of sets~$D$ whose boundaries
contain arbitrarily long arithmetic progressions in the integer lattice~$\Z^2$.

The validity of inequality~\eqref{eq:YudinIntegers}
in the exterior case
is Yudin's theorem; we give a new proof of it
in Section~\ref{sec:visual}.
The fact
that the inequality also holds in the interior case
seems to be new;
we prove it in a
direct
way in Section~\ref{sec:visual},
and
outline a dual proof
in Section~\ref{sec:OtherMethods}.
We
explain in Section~\ref{sec:Convolutions}
how both cases have single-variable precedents in Yves Meyer's paper~\cite{Mey} and related work. We describe the common part of our direct proofs of the two cases
in Section~\ref{sec:TwoLemmas}, and discuss refinements of those methods in Section~\ref{sec:fewer}.
In an appendix, we
outline proofs of two known lemmas that we use
throughout the paper.


The
restriction theorem above
applies to a subspace of~$L^1(\T^2)$ defined by requiring that
some
Fourier coefficients vanish. Related conclusions
hold~\cite{Cooke}, \cite[Theorem~1]{Zyg}
without the latter requirement when~$L^1(\T^2)$ is replaced by~$L^{p}(\T^2)$, where~$p \ge 4/3$.
Unlike most Fourier restriction theorems, that result and ours give global~$\ell^2$ estimates rather than local~$L^2$ estimates.

In Section~\ref{sec:shifted}, we consider
examples
where our methods also yield~$\ell^2$ estimates on  suitable subsets of shifted copies of~$\Gamma$.
These
sometimes lead to global~$L^2$ estimates of the following kind.

\begin{example}
\label{ex:Parabola}
Let~$D = \{(u, v) \in \R^2: v > u^2\}$,
let~$f \in L^1(\R^2)$,
and let~$k \in \R$.
If~$\hat f$ vanishes on~$\Int(D)$ or~$\Int(D^c)$, then
\begin{equation}
\label{eq:BothGlobalL2}
\left[
\int_{-\infty}^\infty |\hat f(u, u^2 + k)|^2\,du
\right]^{1/2}
\le C'|k|^{1/4}\|f\|_1.
\end{equation}
\end{example}


\section
{Contagion of weakness of size
in Fourier algebras
}
\label{sec:Convolutions}

The standard notation for the set of Fourier coefficients of functions in~$L^1(\T^2)$ is~$A(\Z^2)$. This set is
a Banach algebra under pointwise operations
because~$L^1(\T^2)$ is a Banach algebra under convolution.
The
norm of~$\hat f$ in~$A(\Z^2)$ 
is defined to be~$\|f\|_1$. 
Denote the restriction of~$\hat f$ to a set~$S$ by~$\hat f|S$, and
rewrite
inequality~\eqref{eq:YudinIntegers}
in the form
\begin{equation}
\label{eq:FourierRestrict}
\left\|\hat f|(\Gamma\cap\Z^2)\right\|_2
\le C\|\hat f\|_{A(\Z^2)}.
\end{equation}
Also view~$A(\Z^2)$ as the set of sequences on~$\Z^2$ that factor as convolution products of sequences in~$\ell^2(\Z^2)$; this corresponds to the fact that~$L^1 = L^2\cdot L^2$ pointwise.
Moreover,~$\|\hat f\|_{A(\Z^2)}$
is the infimum of the products~$\|g\|_2\|h\|_2$ over all pairs~$(g,h)$ of sequences on~$\Z^2$ for which~$g*h = \hat f$.
Given such a convolution factorization
of~$\hat f$,
extend those factors
the discrete group~$\Rdtwo$ by letting them vanish off~$\Z^2$.
The corresponding extension of~$\hat f$ belongs to~$A(\Rdtwo)$,
with a norm that is
clearly
no larger than the norm of~$\hat f$ in~$A(\Z^2)$.

Theorem~\ref{th:Yudin} follows immediately from the next statement.

\begin{theorem}
\label{th:YudinExtended}
There is a constant~$C$ so that
if~$D$ is
a strictly convex set in~$\R^2$ with boundary~$\Gamma$,
and if~$w \in A(\Rdtwo)$,
then
the estimate
\begin{equation}
\label{eq:Yudin}
\|w|\Gamma\|_2 \le C\|w\|_{A(\Rdtwo)}
\end{equation}
follows from either of the following conditions:
\begin{enumerate}
\item{}
\label{it:Interior}
$w$ vanishes
on~$\Int(D)$.
\item{}
\label{it:Exterior}
$w$ vanishes
on~$\Int(D^c)$.
\end{enumerate}
\end{theorem}

Here we use the notion of ``boundary''
in the usual topology on~$\R^2$.
This makes the corresponding statement for the space~$A(\R^2)$ true but trivial, because functions in~$A(\R^2)$ are continuous
relative to the usual topology on~$\R^2$,
and
they
vanish on~$\Gamma$ if they 
do so
on~$\Int(D)$ or~$\Int(D^c)$.

Meyer's result~\cite[pp. 532--533]{Mey} on~$\Z$
extends
to~$\R_d$
as follows.

\begin{theorem}
\label{th:MeyerExtended}
Let~$(x_j)_{j = 1}^{J}$ be a sequence of positive numbers satisfying the
condition that~$x_{j+1} \ge (1+\delta)x_j$ for some positive constant~$\delta$ and all~$j$.
Let~$w \in A(\Rd)$.
Then
an
estimate
\begin{equation}
\label{eq:Meyer}
\left[\sum_{j = 1}^{J}
|w(x_j)|^2\right]^{1/2}
\le C(\delta)\|w\|_{A(\Rd)}.
\end{equation}
follows from either of the following conditions:
\begin{enumerate}
\item{}
\label{it:before}
$w$ vanishes
on each of the intervals~$\left(
x_j/(1+\delta),
x_j\right)$.
\item{}
\label{it:after}
$w$ vanishes
on each of the intervals~$(x_j, (1+\delta)x_j)$.
\end{enumerate}
\end{theorem}
We will not prove this here, but we  note that,
as in~\cite{FouJMAnAp}, the first part,
about coefficients after
long-enough
gaps, follows by the method that we use to prove the first part of~Theorem~\ref{th:YudinExtended}.
As in~\cite[page~214]{FouArk},
the second part above follows
from
Remark~\ref{rm:RieszSchur} below.

Meyer
used other methods
to prove the version of Theorem~\ref{th:MeyerExtended} for~$A(\Z)$.
He described the pattern in his theorem as a ``contagion of weakness of size.''
On any infinite discrete abelian group~$\Gd$, use the~$\ell^2$ norm
to measure this weakness,
noting that
$
\|w\|_{A(\Gd)} \le \|w\|_2,
$
and recalling that the most one generally say about the size of a function in~$A(\Gd)$ is that it belongs to~$c_0(\Gd)$, which strictly includes~$\ell^2(\Gd)$.

Denote the indicator function of a set~$S$ by~$1_S$.
If~$w|\Int(D)$
belongs to~$\ell^2(D)$,
then
applying
the first part of
Theorem~\ref{th:YudinExtended} to~$w - w\cdot1_{\Int(D)}$
yields that
\begin{equation}
\label{eq:WeakSignalInner}
\|w|\Gamma\|_2
\le C\|w\|_{A(\Rdtwo)}
+ C\|w|\Int(D)\|_2.
\end{equation}
Similarly, if~$w|\Int(D^c)$
belongs to~$\ell^2(D^c)$, then
\begin{equation}
\label{eq:WeakSignalOuter}
\|w|\Gamma\|_2
\le C\|w\|_{A(\Rdtwo)}
+ C\|w|\Int(D^c)\|_2.
\end{equation}
In the setting of Theorem~\ref{th:MeyerExtended},
replace~$\Int(D)$ or~$\Int(D^c)$ with
the union of
long-enough
gaps ending or beginning at the numbers~$x_j$.
In each case, weakness of a member of~$A(\Rdtwo)$ or~$A(\Rd)$ on
a suitable set
propagates to the boundary
of that set
in~$\R^2$ or~$\R$.

\begin{remark}
\label{rm:Paley}
The methods
for the second part of Theorem~\ref{th:YudinExtended}
can also be used~\cite{FouMissing, FouArk, Smith}
to prove Paley's  theorem
about
coefficients of functions
in the classical space~$H^1(\T)$.
In that setting,
weakness on any Hadamard set of positive integers follows from weakness on the set ~$\Z_-$ of negative integers.
It is 
less
clear how
Hadamard
sets in~$Z_+$ can be regarded
as parts of some boundary of~$\Z_-$.
But they share
with the strictly-convex examples
the property
that certain combinations
of  ``boundary points''
must belong to the set
where weakness
is assumed to occur.
See Remark~\ref{rm:whyPaley}
for more on this.
\end{remark}

\begin{remark}
\label{rm:duality}
Recall that~$\Rd$ is dual to the Bohr compactification~$\BohrR$
of the real line.
As in \cite[page 534]{Mey}, applying standard duality arguments  to Theorem~\ref{th:MeyerExtended}
yields that
if~$(v_j) \in \ell^2$, then
there exist
functions~$G$
and~$H$ in~$L^\infty(\BohrR)$ with the following properties.
\begin{enumerate}
\item{}
$\|G\|_\infty$ and~$\|H\|_\infty$ are both no larger than~$C(\delta)\|v\|_2$.
\item{}
$\hat G(x_j) = \hat H(x_j) = v_j$ for all~$j$.
\item{}
$\hat G$ vanishes outside the union
of the intervals~$(x_j/(1+\delta), x_j]$.
\item{}
$\hat H$ vanishes outside the union
of the intervals~$[x_j, (1+\delta)x_j)$.
\end{enumerate}
\end{remark}
\noindent
If~$x_{j+1} \ge (1 + \delta)^2x_j$ for all~$j$, then
the supports of~$\hat G$ and~$\hat H$
are
disjoint except for the numbers~$x_j$.
Work by Goes~\cite[\S 4]{Goes}
exhibited
similar
patterns in a
different
context.
As in~\cite[pp. 214--215]{FouArk},
they
yield an easy proof of
the Grothendieck inequality,
which
follows in the same way from the duals
of
Theorem~\ref{th:Yudin} and~\ref{th:YudinExtended}
that we discuss in Section~\ref{sec:OtherMethods}.

\section{Two Lemmas}
\label{sec:TwoLemmas}

In
our proofs of the
nontrivial
cases
of
Theorem~\ref{th:YudinExtended}, we write each value of~$w$ as an inner product of one function in~$\ell^2(\Rdtwo)$ with a translate of another such function.
Recall that for a
function~$v$ on an additive abelian group and a point~$x$ in that group, the function~$\tau_x v$ maps each point~$y$ to~$v(y-x)$, and the function~$v^*$ maps each point~$y$ to~$\overline{v(-y)}$.
Rename the factor~$h$ in~$w=g*h$ as~$h^*$,
with no effect on norms.
Since
\begin{gather}
(g*h^*)(x)
= \sum_{y\in \R^2} g(y)h^*(x-y)
= \sum_{y\in \R^2} g(y)\overline{h(y-x)},
\notag
\label{eq:gstarhstar}
\\
\label{eq:innerproduct}
w(x) = (g, \tau_x h).
\end{gather}
Proving Theorem~\ref{th:YudinExtended}
therefore
reduces to bounding
$\sum_{j=1}^J  |(g, \tau_{x_j} h)|^2$
for finite sequences~$(x_j)_{j=1}^J$ of distinct  points in~$\Gamma$.

We apply the lemmas below with~$H = \ell^2(\Rdtwo)$ and~$A_j = \tau_{x_j}$. The first lemma goes back to~\cite{FouPac}, and led to a rediscovery~\cite{FouJMAnAp} of Meyer's result about coefficients after gaps. The second lemma is more recent~\cite{FouMissing},
and was used
there
to reprove 
the
extension~\cite[Theorem~2]{FouArk} of
Paley's theorem
that
yields the part of Theorem~\ref{th:MeyerExtended}
about coefficients before gaps. In the next section, we specify subspaces with the properties required in the lemmas. We outline proofs of the lemmas in Appendix~\ref{sec:TwoStep}.

\begin{lemma}\label{th:oldHilbert}
Let $H$ be a Hilbert space
and~$M_1 \subset M_2 \subset \cdots \subset M_J$
be closed subspaces of~$H$.
Let~$A_1, A_2, \dots, A_J$ be unitary operators on~$H$
for which
\[
A_1M_1 \subset A_2M_2 \subset \cdots \subset A_JM_J.
\]
Let~$g$ and~$h$ be
members
of~$H$
satisfying the following conditions
for all indices~$j < J$.
\begin{enumerate}
\item\label{itF:intarget}
$A_jh \in A_{j+1}M_{j+1}$.
\item\label{itF:orthogonal}
The vector $g$ is orthogonal to the subspace~$A_{j+1}M_j$.
\end{enumerate}
Then
\begin{equation}
\label{eq:oldHilbertBound}
\left[\sum_{j=1}^J |(g, A_jh)|^2\right]^{1/2}
\le 2(\|g\|_H)
\|h\|_H.
\end{equation}
\end{lemma}

\begin{lemma}\label{th:newHilbert}
Let $H$ be a Hilbert space
and
$
L_1 \supset L_2 \supset \cdots \supset L_J
$
be closed subspaces of~$H$.
Let~$A_1, A_2, \dots, A_J$ be unitary operators on~$H$
for which
\[
A_2L_1 \subset A_3L_2 \subset \cdots \subset A_JL_{J-1}.
\]
Let~$g$ and~$h$ be elements of~$H$
satisfying the following conditions:
\begin{enumerate}
\item\label{it:membership}
$A_jh \in A_{j+1}L_j$ for all~$j < J$.
\item\label{it:orthogonality}
The vector~$g$ is orthogonal to the subspace~$A_{j}L_j$ for all~$j > 1$.
\end{enumerate}
Then
\begin{equation}
\label{eq:newHilbertBound}
\left[\sum_{j=1}^J |(g, A_jh)|^2\right]^{1/2}
\le 2(\|g\|_H)
\|h\|_H.
\end{equation}
\end{lemma}

\section{Visual proofs}
\label{sec:visual}


Given
the convolution factorization~$w = g*h^*$
and
a subset~$E$
of~$\Rdtwo$, let~$V(E, h)$ denote the closure in~$H = \ell^2(\Rdtwo)$ of the subspace spanned by the translates~$\tau_x h$ for which~$x \in E$.
In the interior case
of Theorem~\ref{th:YudinExtended}, we
will
apply
Lemma~\ref{th:oldHilbert}
with~$M_j = V(E_j, h)$
for suitable sets~$E_j$.
In the exterior case, we
will
apply Lemma~\ref{th:newHilbert} with~$L_j = V(D_j, h)$ for suitable sets~$D_j$.

The nesting and membership conditions in
Lemma~\ref{th:oldHilbert}
hold if
\begin{equation}
E_j \subset E_{j+1},
\quad
x_j + E_j \subset x_{j+1} + E_{j+1},
\quad\text{and~$x_j \in x_{j+1} + E_{j+1}$}
\end{equation}
for all~$j < J$.
The
orthogonality condition
holds if~$(g, \tau_y h) = 0$ for all~$y$ in~$x_{j+1} + E_j$.
Equation~\eqref{eq:innerproduct}
makes this
equivalent to
having~$w(y) = 0$
for all such~$y$.

Let~$F_j = x_j + E_j$ for all~$j$,
and let~$\Delta x_j = x_{j+1} - x_j$ when~$j < J$.
The
last condition in the previous paragraph
is
equivalent to requiring
that~$w$ vanish on
all
the sets~$F_j + \Delta x_j$ with~$j < J$.
In
the interior case,
this happens if
those sets are all included in~$\Int(D)$.
Translate the other conditions on the sets~$E_j$
to see that it suffices
in that case
to find
sets~$F_j$ satisfying the following four conditions
for all~$j < J$.
\begin{gather}
F_j + \Delta x_j \subset \Int(D),
\quad
F_j + \Delta x_j \subset F_{j+1},
\label{eq:Fj+Deltaxj}
\\
F_j \subset F_{j+1},
\quad\text{and}\quad
x_j \in F_{j+1}.
\label{eq:F_j+1}
\end{gather}
That is,
\begin{equation}
\label{eq:InnerInclusions}
F_j + \Delta x_j \subset \Int(D)\cap F_{j+1}
\quad\textrm{and}\quad
F_j\cup\{x_j\} \subset F_{j+1}.
\end{equation}
Call these the \textit{shifted inclusions} and the \textit{unshifted inclusions}.

Similarly, the subspaces~$L_j$ and their images~$A_{j+1}L_j$ nest
as prescribed in Lemma~\ref{th:newHilbert}
if
\begin{gather}
\label{eq:antinest}
D_1 \supset D_2 \supset \cdots \supset D_J,
\notag
\\
\label{eq:imagenest}
\text{and}\quad
x_2 + D_1 \subset x_3 + D_2 \subset \cdots
\subset x_J + D_{J-1}.
\notag
\end{gather}
The membership
condition in the lemma
holds if
$
x_j \in x_{j+1} + D_j
$
for all~$j < J$,
and
the orthogonality condition holds
in the exterior case
if~$
x_j + D_j \subset \Int(D^c)
$
for all~$j > 1$.

Consider the
sets~$G_{j+1} = x_{j+1}+ D_j$,
creating another point~$x_{J+1}$
to cover the case where~$j = J$.
Translate the
conditions on the sets~$D_j$
to see that it suffices
that
\begin{gather}
\label{eq:Gj-Deltaxj}
G_{j+1} - \Delta x_j \subset \Int(D^c),
\quad
G_{j+1} - \Delta x_j \subset G_{j},
\\
\label{eq:G_j+1}
G_{j} \subset G_{j+1},
\quad\text{and}\quad
x_j \in G_{j+1}.
\end{gather}
In this case,
the shifted inclusions and unshifted inclusions state that
\begin{equation}
\label{eq:OuterInclusions}
G_{j+1} - \Delta x_j \subset \Int(D^c)\cap G_{j}
\quad\textrm{and}\quad
G_j\cup\{x_j\} \subset G_{j+1}.
\end{equation}

If the boundary~$\Gamma$ of~$D$ is the graph of a
strictly convex or
strictly concave function defined on all of~$\R$,
and the points~$x_j$ run from left to right along~$\Gamma$,
then we can use sets~$F_j$ and~$G_j$ that are very similar.
For such a concave function~$\phi$,
write~$x_j = (u_j, v_j)$, and
\begin{gather}
\textrm{let}\quad
F_j = \{(u,v) \in \R^2: u < u_j
\textrm{ and }
v \le \phi(u)\},
\label{eq:SimilarFj}
\\
\textrm{while}\quad
G_j = \{(u,v) \in \R^2: u < u_j
\textrm{ and }
v \ge \phi(u)\}.
\label{SimilarGj}
\end{gather}
\begin{figure}[h]
        \centering
        \begin{subfigure}[b]{0.5\textwidth}
                \centering
                \includegraphics[width=145 pt]{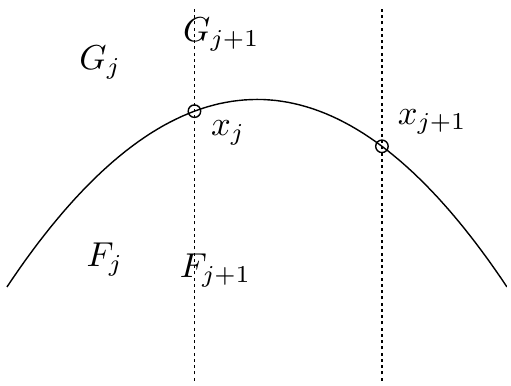}
                \caption{Unshifted
                Inclusions
                }
                \label{fi:BothUnshifted}
        \end{subfigure}%
        \qquad
        ~ 
        \begin{subfigure}[b]{0.5\textwidth}
                \centering
                \includegraphics[width=145 pt]{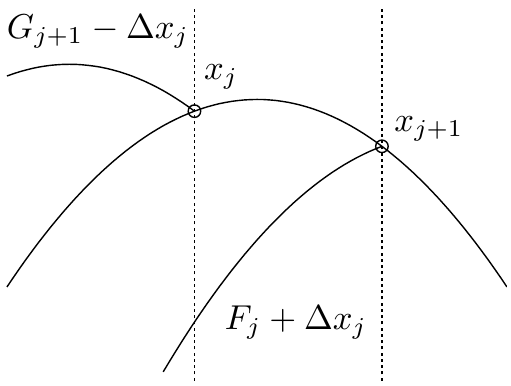}
                \caption{Shifted
               Inclusions
                }
                \label{fi:BothShifted}
        \end{subfigure}
\caption{Similar Sets}
\label{Fi:Both}
\end{figure}

\noindent
The unshifted inclusions
clearly
hold for
both~$F_j$ and~$G_j$.
By
strict concavity, any part of the
boundary ending at~$x_{j}$ rises strictly more rapidy or falls strictly more slowly than any part of the same width to the right of it. Shifting such a part
ending at~$x_j$
by~$\Delta x_j$
gives
a curve that ends at~$x_{j+1}$ and lies strictly below~$\Gamma$ except at~$x_{j+1}$. This
yields
the shifted inclusions for the sets~$F_j + \Delta x_j$. The corresponding inclusions for the sets~$G_{j+1} - \Delta x_j$ follow in
a similar way.
Both
cases of Theorem~\ref{th:YudinExtended} 
therefore
hold with~$C = 2$
for such sets~$D$.

Every unbounded, strictly convex
set
can be rotated to have the form specified above, except that
the domain of the function~$\phi$
may not be
all of~$\R$.
In that case,
add the requirement that~$u$ belong to the domain of~$\phi$ in defining~$F_j$. If the domain
of~$\phi$ 
is bounded 
on the left, 
also include all vertical lines to the left 
of~$D$
in defining~$G_j$.

When~$D$
is bounded and strictly convex,
follow~\cite{Yudin}
in recalling that
there are vertical support lines
at two boundary points, listed from left to right as~$x_0$ and~$x_\infty$ say.
In the exterior case,
let~$\Gamma_0$ be the upper boundary
with~$x_\infty$ excluded.

Consider points~$x_j$ running from left to right in~$\Gamma_0$, starting with~$x_0$.
As above,
let~$G_j$
consist of all points in~$\R^2$ that lie strictly to the left of~$x_j$, and that do not lie directly below~$\Gamma_0$.
Then the 
inclusions~\eqref{eq:OuterInclusions} 
hold 
for all~$j \ge 0$,
so that~$\left\|w|\Gamma_0\right\|_2
\le 2 \|w\|_{A(\Rdtwo)}$.
Rotate by~$180^\circ$ to get
a similar 
estimate on the rest 
of~$\Gamma$,
and
that
\begin{equation}
\label{eq:VanishOutside}
\|w|\Gamma\|_2 \le
2\sqrt2\|w\|_{A(\Rdtwo)}.
\end{equation}

In the interior case
for the same set~$D$,
shear vertically and shift
to place both of
the
points~$x_0$ and~$x_\infty$
on the~$u$-axis;
this does not affect~$\|w\|_{A(\Rdtwo)}$.
Then the lower boundary lies below the~$u$-axis.
There will be one point,~$x_{J}$ say,
on the upper boundary
with a horizontal support line.
Place that point
on the~$v$-axis.
%
Then the upper boundary in the second quadrant is
the graph of
an
increasing function.

Consider points~$\{x_j\}_{j=1}^{J-1}$
running from left to right
in the interior of
that graph.
Find the midpoint of the  line segment from~$x_0$ to~$x_j$; then rotate the part of boundary curve running from~$x_0$ to~$x_j$ by~$180^o$ about that midpoint
to get
a lower curve
returning
to~$x_0$
from~$x_j$.
Form the convex hull of that lower curve and the upper boundary curve from~$x_0$ to~$x_j$,
and
delete the
vertices~$x_0$ and~$x_j$ to get the set~$F_j$.
Form~$F_J$ in the same way.
We show this in
Figure~\ref{fi:InnerUnshifted}
below.
\begin{figure}[h]
        \centering
        \begin{subfigure}[b]{0.5\textwidth}
                \centering
                \includegraphics[width=145 pt]{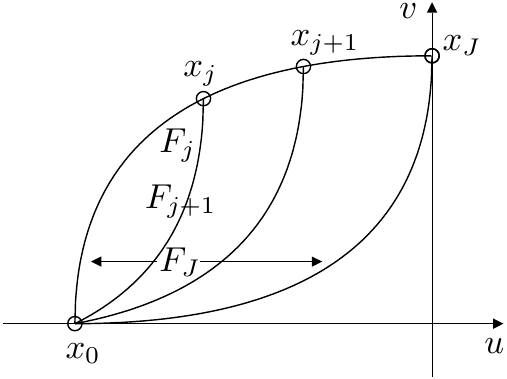}
                \caption{$\{x_{j}\}\cup F_{j}\subset F_{j+1}$}
                \label{fi:InnerUnshifted}
        \end{subfigure}%
        \qquad
        ~ 
        \begin{subfigure}[b]{0.5\textwidth}
                \centering
                \includegraphics[width=145 pt]{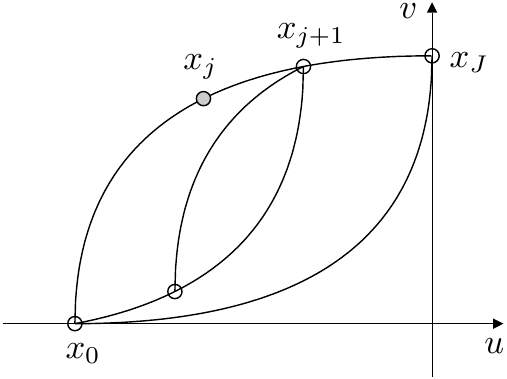}
                \caption{$
                F_j + \Delta x_j 
                \subset \Int(D) \cap F_{j+1}$}
                \label{fi:InnerShifted}
        \end{subfigure}
\caption{Interior Inclusions}
 \label{Fi:Inner}
\end{figure}

\noindent
It is obvious that~$\{x_j\} \subset F_{j+1}$.
The sets~$F_j$ and~$F_{j+1}$ are mapped onto themselves by the~$180^\circ$ rotations,~$\psi_j$ and~$\psi_{j+1}$ say, about their centroids.
Note that~$\psi_{j+1}$ is equal to~$\psi_j$ followed by the shift by~$\Delta x_j$;
so~$\psi_{j+1}$
maps~$F_j$ onto~$F_j + \Delta x_j$.
Since the upper boundary of~$F_j$ 
is, by definition,
an initial part of the upper boundary of~$F_{j+1}$,
the lower boundary of~$F_j + \Delta x_j
= \psi_{j+1}(F_j)$
is
a final part of the lower boundary of~$F_{j+1}$.

As in Figure~\ref{fi:BothShifted},
the upper boundary
of~$F_j + \Delta x_j$ lies strictly below the upper boundary of~$F_{j+1}$ except at the 
missing 
point~$x_{j+1}
$.
Hence~$F_j + \Delta x_j \subset F_{j+1}$; applying~$\psi_{j+1}$ again then makes~$F_j  \subset F_{j+1}$.

The
lower boundary
of~$F_J$
runs from~$x_0$ to~$x_J$,
and is
the graph of
an
increasing
function.
Hence~$F_J$ lies strictly inside the second quadrant, as do its subsets~$F_j + \Delta x_j$ with~$j < J$. 
These
sets
therefore
do not meet
the lower boundary 
or right-hand boundary
of~$D$. 
Since 
the shifted sets~$F_j + \Delta x_j$
lie
strictly
below
the
upper boundary
of~$D$
in the second quadrant,
they 
are included in~$\Int(D)$,
as required.

Let~$\Gamma_2$ be the part of~$\Gamma$ inside the second quadrant, together with~$x_J$.
Then~$\left\|w|\Gamma_2\right\|_2
\le 2 \|w\|_{A(\Rdtwo)}$ in the interior case.
Similar
arguments
on three other parts of~$\Gamma$
yield
that
\begin{equation}
\label{eq:VanishInside}
\|w|\Gamma\|_2 \le
4\|w\|_{A(\Rdtwo)}.
\end{equation}

\section{Weaker hypotheses}
\label{sec:fewer}

Our
methods work when~$w$ vanishes on
some sets that are smaller than the ones used
in Section~\ref{sec:visual}.
In the next section, we 
discuss
dual methods that also work with those weaker hypotheses.

Fix a finite sequence~$(x_j)_{j=1}^J$. It will turn out to suffice that~$w$ vanish on suitable subsets of the additive group generated by the points~$x_j$.
All points~$x$ in that group have the form
\begin{equation}\label{eq:epsilons}
x = \sum_{i=1}^{J} \varepsilon_i x_i,
\end{equation}
where the coefficients~$\varepsilon_i$ are integers.

The application of Lemma~\ref{th:oldHilbert}
to lacunary Fourier series was analysed in~\cite[Remark~3]{FouJMAnAp}.
In the present context, the same reasoning shows that it suffices for~$w$ to vanish on the set~$\Alt((x_j))$ of points~$x$ with alternating sum representations
\begin{equation}
\label{eq:Alt}
x =
x_{j_1} - x_{j_2} + \cdots +
x_{j_{2i-1}} - x_{j_{2i}} + x_{j_{2i+1}}
\end{equation}
with at least~$3$ terms and a strictly-increasing index sequence~$(j_\ell)$.

Let~$F_{j+1}$ be the set of points~$x$ as above with~$j_{2i+1} \le j+1$, but only impose the requirement that the sum~\eqref{eq:Alt} have at least~$3$ terms when~$j_{2i+1} = j+1$.
These sums belong to the fatter sets~$F_{j+1}$ shown
in
Figures~\ref{fi:BothUnshifted} and~\ref{fi:InnerUnshifted}.
The inclusions~\eqref{eq:InnerInclusions}
hold for the smaller sets~$F_j$ and~$F_{j+1}$, and Lemma~\ref{th:oldHilbert} applies.

For Lemma~\ref{th:newHilbert}, the analysis in~\cite[Section~5]{FouMissing}
yields the sets~$G_{j+1}$ consisting of
all points~$x$
with
a representation
\begin{equation}
\label{eq:DeltaXj's}
x = x_i - \sum_{j' \ge i} n_{j'}\Delta x_{j'},
\end{equation}
satisfying the following conditions:
\begin{enumerate}
\item{}
\label{StartLow}
$i \le j+1$.
\item{}
\label{NnIntegers}
The coefficients~$n_{j'}$ are nonnegative integers.
\item{}
\label{Consecutive}
If~$i = j+1$, then~$n_{j'} \ne 0$ for some~$j'$.
\end{enumerate}
The points in this version of~$G_{j+1}$ belong to
the fatter set~$G_{j+1}$ shown in
Figure~\ref{fi:BothUnshifted}.
The
desired
inclusions
hold
for the smaller sets~$G_j$ and~$G_{j+1}$.

The lemma applies provided that~$w$ vanishes on the union~$\Schur((x_j))$ of the smaller sets~${G_{j+1}} - \Delta x_j$. The points~$x$ in that union are those with a representation~\eqref{eq:DeltaXj's} satisfying condition~\eqref{NnIntegers} with~$n_{j'} \ne 0$ for some~$j'$.
They are also
given by
the
sums of the form~\eqref{eq:epsilons}
where the
integer
coefficients~$\varepsilon_i$ have the following
properties:
\begin{itemize}
\item{}
\label{it:fullsum}
The full sum~$\sum_{i=1}^{J} \varepsilon_i $
is equal to~$1$.
\item{}
All partial sums of the full sum are nonnegative.
\item{}
All partial sums after the first positive one are positive.
\item{}
Some partial sum is greater than~$1$.
\end{itemize}
These conditions
also arose
in~\cite{FouArk} and~\cite{Yudin}.

\begin{remark}
\label{rm:whyPaley}
In the setting of Remark~\ref{rm:Paley},
Paley's theorem holds because the set~$\Schur(\{n_j\})$
is included
in~$\Z_-$
when the sequence~$(n_j)$ is sufficiently lacunary.
This was used in a dual way in~\cite{FouArk}
and~\cite{Smith},
and in a direct way in~\cite{FouMissing}.
\end{remark}

\begin{remark}
\label{rm:PartialOrders}
We made one choice of the sets~$G_j$ in proving the exterior case of Theorem~\ref{th:YudinExtended},
and
another
just above. For
both
choices, the corresponding sets~$D_j$ are additive semigroups. This can be used~\cite[Remark~5.7]{FouMissing} to define suitable partial orders 
on~$\R^2$,
relating
that case of Theorem~\ref{th:YudinExtended} to Paley's theorem.
\end{remark}

\begin{remark}
\label{rm:RieszSchur}
It can happen that~$|\varepsilon_i| > 1$ in the sums~\eqref{eq:epsilons} representing points in~$\Schur((x_j))$. Let~$S((x_j))$ consist of all points~$x$ with representations~\eqref{eq:epsilons} in which the coefficients~$\varepsilon_{i}$ belong to the set~$\{-1, 0, 1\}$ and satisfy the four conditions for membership of~$x$ in~$\Schur((x_j))$. Arguments in~\cite{FouArk} and~\cite{FouMissing} each combine with the application above of Lemma~\ref{th:newHilbert} to show that
\begin{equation}
\|w|X\|_2 \le 4\|w\|_{A(\Rdtwo)}
\end{equation}
when~$w$ vanishes on~$S((x_j))$.
\end{remark}

\begin{remark}
\label{rm:DeltaFj}
For~$\Alt((x_j))$,
rewrite
the representation~\eqref{eq:Alt}
in the form
\begin{equation}
\label{eq:AltDeltaXj's}
x= x_{j_{2i+1} }
- \sum_{j' < j_{2i+1} - 1}
n_{j'}\Delta x_{j'}
\end{equation}
where the coefficients~$n_{j'}$ take the values~$0$ and~$1$ only
and the latter occurs at least once.
For~$F_{j+1}$, keep
those
conditions
on~$(n_{j'})$,
put~$j_{2i+1} = j+1$, and
require
instead
that~$j' \le j$
in the sum.
\end{remark}

\section{Dual
Constructions
}
\label{sec:OtherMethods}

Denote the Bohr compactification of~$\R^2$ by~$\BohrR^2$.
The duality arguments in~\cite{RudPal} or~\cite[page 534]{Mey} show that
Theorem~\ref{th:YudinExtended}
is equivalent to the
one below.
Theorem~\ref{th:Yudin}
has a similar dual.

\begin{theorem}
\label{th:DualYudin}
Let~$D$ be a strictly convex set in~$\R^2$
with boundary~$\Gamma$.
Then
for each function~$v$ in~$\ell^2(\Gamma)$, there exist functions~$G$
and~$H$ in~$L^\infty(\BohrR^2)$ with the following properties:
\begin{enumerate}
\item{}
$\|G\|_\infty$ and~$\|H\|_\infty$ are both no larger than~$C\|v\|_2$.
\item{}
The restrictions of~$\hat G$ and~$\hat H$ to~$\Gamma$  both coincide with~$v$.
\item{}
$\hat G$ vanishes on~$\Int(D^c)$.
\item{}
$\hat H$ vanishes on~$\Int(D)$.
\end{enumerate}
\end{theorem}

Theorem~\ref{th:YudinExtended}
can
be proved by constructing suitable functions~$G$ and~$H$ when 
the support of~$v$ is finite. 
Choose points~$x_j$ as in Section~\ref{sec:visual}. Let~$v$ vanish off the set~$X = \{x_j\}_{j=1}^J$, with~$\|v\|_2 =1$.
The
modification of the Rudin-Shapiro construction in~\cite{Clunie} 
produces a trigonometric polynomial~$G$ with the following properties.
\begin{itemize}
\item{}
$\|G\|_\infty \le C$.
\item{}
$\hat G|X = v$ if
the sets~$X$ and~$\Alt((x_j))$ are disjoint.
\item{}
$\hat G$ vanishes off the set~$X\cup \Alt((x_j))$.
\end{itemize}
This yields the first part of Theorem~\ref{th:YudinExtended}, since the strict convexity of the unbounded set~$D$ makes~$\Alt((x_j))$ a subset of~$\Int(D)$
in the diagrams in Section~\ref{sec:visual}.

It also follows that~$\Schur((x_j)) \subset \Int(D^c)$
in those cases.
For the second
part of the theorem,
it suffices to
construct a function~$H$ with
the following properties.
\begin{itemize}
\item{}
$\|H\|_\infty \le 1$.
\item{}
$\overline{v(x_j)}{\hat H(x_j)} \ge
(1/C)|v(x_j)|^2$
for all~$x_j$.
\item{}
$\hat H$ vanishes off the set~$X\cup \Schur((x_j))$.
\end{itemize}
Yudin refined a method of Pigno and Smith~\cite{PS, Smith}
for
this, and
noted that a
construction
in~\cite{FouArk}
would work too.
In
both
of
these methods,
one can satisfy
the middle condition above
by making~$\hat H|K = (1/C)v$.

\section{Separated points in shifted curves}
\label{sec:shifted}

In
Example~\ref{ex:Parabola},
let~$k > 0$ and~$h
>
\sqrt{k/2}$.
We will show
that
\begin{equation}
\label{eq:InnerAmalgam}
\left[
\sum_{j=-\infty}^\infty \left\{
\sup_{jh \le u < (j+1)h}|\hat f(u, u^2 - k)|
\right\}^2
\right]^{1/2}
\le \sqrt2C\|f\|_1,
\end{equation}
in the interior case, and that
\begin{equation}
\label{eq:OuterAmalgam}
\left[\sum_{j=-\infty}^\infty \left\{\sup_{jh \le u < (j+1)h}|\hat f(u, u^2 + k)|\right\}^2\right]^{1/2}\le \sqrt2C\|f\|_1
\end{equation}
in the
exterior case.
Inequality~\eqref{eq:BothGlobalL2} 
then
follows
because
\[
\int_{jh}^{(j+1)h} |g(u)|^2 \,du
\le h\left\{\sup_{jh \le u < (j+1)h} |g(u)|\right\}^2
\]
for all measurable functions~$g$.

The ``amalgam norm'' estimates~\eqref{eq:InnerAmalgam} and~\eqref{eq:OuterAmalgam}
follow from~$\ell^2$ estimates 
on sets of suitably separated points,~$x_j = (u_j, u_j^2)$ say,
in~$\Gamma$.
Let~$w \in A(\R^2)$
and
require 
it to vanish
that~$w$ vanishes
on the region where~$v > u^2 + k$,
or on the region where~$v < u^2 - k$.
Then
\begin{equation}
\label{eq:ShiftedEstimate}
\left[
\sum_j |
w(x_j)|^2
\right]^{1/2}
\le C\|w\|_{A(\R^2)}
\quad\textnormal{if~$\Delta u_j > \sqrt{k/2}$ for all~$j$.}
\end{equation}
Apply this to shifted copies~$w$ of~$\hat f$, and choose points~$x_j$ in alternate intervals~$[j'h, (j'+1)h)$ to get
the estimates~\eqref{eq:InnerAmalgam} and~\eqref{eq:OuterAmalgam}.

In proving inequality~\eqref{eq:ShiftedEstimate},
we consider more general sets~$D$ of the form~$\{(u,v): v \ge \phi(u)\}$,
where~$\phi''(u) \ge 
c 
> 0$.
Our methods apply 
to~$A(\R^2)$, and
yield inequality~\eqref{eq:ShiftedEstimate} if the sets~$\Alt((x_j))$ and~$\Schur((x_j))$ are respectively included in the
sets~$\Int(D)+(0,k)$ and~$\Int(D^c)-(0,k)$.

Given
a point~$x$ in~$\Schur((x_j))$
in
the form~\eqref{eq:DeltaXj's},
let~$n = \sum_{j'}n_{j'}$
and say
that~$x$ is an~$n$-th generation descendant of~$x_i$.
Subtracting another copy of~$\Delta x_j$, where~$j \ge i$, from~$x$ gives an~$(n+1)$-st descendant,~$x'$ say.
All
descendants~$(u, v)$
of~$x_i$
share the property that~$u < u_i$.
Visual arguments in the style of Section~\ref{sec:visual} show that
if~$\phi$ is stricly convex
and~$x \in \Int(D^c - (0, k))$,
then~$x' \in \Int(D^c - (0, k))$ too.

So
it suffices
to
check
that 
first-generation points in~$\Schur((x_j))$ belong to~$\Int(D^c-(0,k))$.
They
have the form~$x_i - \Delta x_{j}$ where~$j\ge i$.
Rewriting this as~$(u,v)
= (u_i, v_i) - (\Delta u_j, \Delta v_j)$
reduces matters to
showing
that~$\phi(u) - v > k$.
Now
\begin{gather}
\notag
v_i =
\phi(u_i) = \phi(u) + \int_{u}^{u_i}
\phi'(r)\,dr,
\quad
\phi(u) =
v_i - \int_{0}^{\Delta u_j}
\phi'(u + s)\,ds,
\\
\notag
\Delta v_j = \int_{u_j}^{u_{j+1}}
\phi'(r)\,dr,
\quad\text{and}\quad
v = v_i - \Delta v_j
= v_i - \int_{0}^{\Delta u_{j}}
\phi'(u_j + s)\,ds.
\end{gather}
Therefore,
\begin{gather}
\phi(u) - v =
\int_{0}^{\Delta u_{j}}
\left[\phi'(u_j + s) - \phi'(u + s)\right]
\,ds
\notag
\\
=
\int_{0}^{\Delta u_{j}}
\left[\int_u^{u_j}
\phi''(t+s)
\,dt\right]
\,ds
\ge c(\Delta u_j)^2.
\notag
\end{gather}

Use the representation~\eqref{eq:AltDeltaXj's}
to introduce a similar notion of generations of descendants in~$\Alt((x_j))$,
but add the requirement that the extra nonzero coefficient~$n_{j'}$ for the child~$x' = (u', v')$ occurs before all nonzero coefficients for the parent~$x$.
Rename~$j_{2i+1}$ as~$j+1$;
then~$u' \ge u_{j'} + \Delta u_j$.
Argue visually to
reduce matters to first-generation cases
where
\[
x' = (u', v') =  (u_{j+1}, v_{j+1}) - (\Delta u_{j'}, \Delta v_{j'}),  \quad{\rm and\ }j' < j.
\]
As above,
\[
v' - \phi(u')
=
\int_{0}^{\Delta u_{j'}}
\left[\int_{u_{j'}}^{u'}
\phi''(t+s)
\,dt\right]
\,ds
\ge
c(\Delta u_{j'})\Delta u_{j}.
\]
The
inclusions~$\Schur((x_j)) \subset \Int(D^c - (0,k))$ and~$\Alt((x_j)) \subset \Int(D + (0,k))$
follow
if~$\Delta u_j > \sqrt{k/c}$ for all~$j$.

The outcome changes if the graph of~$\phi$ has an asymptote.

\begin{example}
\label{ex:Hyperbola}
Let~$D_\alpha = \{(u, v): u>0, v > u^{-\alpha}\}$,
where~$\alpha$ is a positive constant.
Let~$f \in L^1(\R^2)$, and let~$k > 0$.
If~$\hat f$ vanishes on~$\Int(D_\alpha^c)$, then
\begin{equation}
\label{eq:WeightedL2}
\left[\int_0^\infty
\left|\hat f\left(u, u^{-\alpha} + k\right)\right|^2
\,\frac{du}{u}\right]^{1/2}
\le
C\|f\|_1.
\end{equation}
There
are cases where~$f \in L^1(\R^2)$
and~$\hat f$ vanishes on~$\Int(D_\alpha)$
but
\[
\int_{0}^\infty
\left|\hat f\left(u, u^{-\alpha} - k\right)\right|^2
\,\frac{du}{u}
= \infty.
\]
\end{example}

The positive result here follows from
the
extension of Paley's inequality to functions~$f$ in~$L^1(\R^2)$ for which~$\hat f(u, v) = 0$ on the ``negative'' semigroup,~$-P$ say, where~$u \le 0$ and~$v < 0$ if~$u = 0$. 
That extension
gives an~$\ell^2$ estimate for~$(\hat f(x_j))$ when the sequence~$(x_j)$ satisfies the Hadamard condition that~$2 x_j - x_{j+1}\in -P$ for all~$j$.
So do the appropriate methods in Sections~\ref{sec:fewer} or~\ref{sec:OtherMethods}.
These approaches all show that
\[
\left[\int_0^\infty
\sup_{v \in \R}
\left|\hat f\left(u, v\right)\right|^2
\,\frac{du}{u}\right]^{1/2}
\le
C\|f\|_1.
\]

To get the negative results, use the fact that for each parallelogram,~$B$ say, with positive area, there is a function in the unit ball of~$A(\R^2)$ that vanishes outside~$B$ and that exceeds~$1/4$ on~$1/4$ of the area of~$B$. One way to confirm this fact runs via the
argument applied to arithmetic progressions in~\cite[p. 861]{Yudin}.

Similar reasoning, going back to~\cite{RudPal}, shows that if a nonnegative measure~$\nu$ has the property that
\[
\int_{0}^\infty
\left|\hat f\left(u, u^{-\alpha} + k\right)\right|^2
\,
d\nu(u)
< \infty
\]
whenever~$f \in L^1(\R^2)$ and~$\hat f$ vanishes on~$\Int(D_\alpha^c)$, then
\[
\nu((2^j, 2^{j+1}]) \le C'
\quad\text{for all~$j$.}
\]

\begin{remark}
\label{rm:Affine}
Affine arclength measure is prominent in restriction theorems~\cite{Ob} for transforms of functions in~~$L^p(\R^2)$ when~$p > 1$.
The measure~$du$
on the graphs of~$v = u^2 \pm k$
is affine invariant, but the measure~$du/u$
on the graph of~$v = \phi_\alpha(u) + k$
is not, except when~$\alpha = 1$.
\end{remark}

\appendix
\section{
Two orthogonality steps}
\label{sec:TwoStep}

We prove both
lemmas 
by splitting
the sequence~$(g, A_jh)_{j=1}^J$ as a sum of two sequences whose~$\ell^2$ norms are easy to bound.

In
Lemma~\ref{th:newHilbert}, let~$P_j$ and~$Q_j$
be the orthogonal projections
onto the subspaces~$L_{j}$ and~$A_{j+1}L_{j}$ respectively,
with~$j < J$ in the latter case. Also let~$Q_J = I$
and~$Q_0 = 0$.
By
the
membership condition in the lemma,
\begin{equation}\label{eq:definition}
(g, A_jh) = (g, Q_{j}A_jh) = (Q_{j}g, A_jh)
= a_j + b_j,
\end{equation}
where~$a_j = ((Q_j - Q_{j-1})g,A_jh)$ and~$b_j = (Q_{j-1}g, A_jh)$ for all~$j$.
Then~$b_1 = (Q_0g, A_jh) = 0$,
and~$b_j = (g,A_j(P_{j-1} - P_j)h)$ when~$j > 1$,
since~$A_jP_{j-1} = Q_{j-1}A_j$ and~$(g, A_jP_jh) = 0$
in 
that
case.
The projections~$Q_{j} - Q_{j-1}$
have mutually orthogonal ranges,
as do the projections~$P_{j-1} - P_{j}$. By Cauchy-Schwarz,~$\|(a_j)||_2$
and~$\|(b_j)\|_2$ are both bounded above by~$(\|g\|_H)\|h\|_H$, and inequality~\eqref{eq:newHilbertBound} follows.

In
Lemma~\ref{th:oldHilbert},
consider the orthogonal projections~$P_j$ and~$Q_j$ onto the subspaces~$M_j$ and~$A_jM_j$.
Also let~$Q_{J+1} = I$ and~$P_0 = 0$.
This time,~$(g, A_jh) = (Q_{j+1}g, A_jh)$,
which
splits as
\begin{equation}
((Q_{j+1}-Q_{j})g,A_jh)
+ (g,A_j(P_j - P_{j-1})h),
\end{equation}
since~$(g, A_jP_{j-1}h) = 0$ and~$(g, A_jP_{j}h) = (Q_{j}g,A_jh)$.
Finish
as above.
\bibliographystyle{amsplain}

\begin{thebibliography}{10}

\bibitem{Clunie}
J. Clunie,
\textit{On the derivative of a bounded function,}
Proc. London Math. Soc.~(3) \textbf{14a} (1965), 58--68.

\bibitem{Cooke}
Roger Cooke,
\textit{A Cantor-Lebesgue theorem in two dimensions},
Proc. Amer. Math. Soc.
\textbf{30} (1971), 547--550.

\bibitem{FouPac}
John J.F. Fournier,
\textit{Extensions of a Fourier multiplier theorem
of Paley,}
Pacific J. Math.
\textbf{30} (1969), 415--431.

\bibitem{FouJMAnAp}
\bysame,
\textit{Fourier coefficients after gaps,}
J. Math. Anal. Appl.
\textbf{42} (1973), 255--270.


\bibitem{FouArk}
\bysame,
\textit{On a theorem of Paley and the Littlewood conjecture,}
Ark. Mat.
\textbf{17} (1979), 199--216.

\bibitem{FouMissing}
\bysame,
\textit{The missing proof of Paley's theorem about lacunary coefficients,}
arXiv:1407.1458 [math.CA].

\bibitem{Goes}
Gunther Goes,
\textit{On a Tauberian theorem for sequences with gaps and on Fourier series with gaps,}
T\^{o}hoku Math. J. \textbf{24} (1972), 153--165.

\bibitem{Mey}
Yves Meyer,
\textit{Endomorphismes des id\'eaux ferm\'es de $L^{1}(G)$, classes de Hardy et s\'eries de Fourier lacunaires,} Ann. Sci. \'Ecole Norm. Sup. (4) \textbf{1} (1968), 499--580.

\bibitem{Ob}
Daniel M. Oberlin,
\textit{Fourier restriction for affine arclength measures in the plane,}
Proc. Amer. Math. Soc. \textbf{129} (2001),
3303--3305.

\bibitem{PS}
Louis Pigno and  Brent Smith,
\textit{A Littlewood-Paley inequality for analytic measures.}
Ark. Mat. \textbf{20} (1982),
271--274.

\bibitem{REACH1}
R.E.A.C. Paley,
\textit{On the lacunary coefficients of power series,}
Ann. of Math.~(2) \textbf{34} (1933), 615--616.

\bibitem{RudPal}
Walter Rudin,
\textit{Remarks on a theorem of Paley,}
J. London Math. Soc.
\textbf{32} (1957), 307--311.

\bibitem{Smith}
Brent Smith,
\textit{Two trigonometric designs: one-sided Riesz products and Littlewood products,}
General inequalities, \textbf{3} (Oberwolfach, 1981), 141--148, Internat. Schriftenreihe Numer. Math., 64, Birkh�user, Basel, 1983. 

\bibitem{Yudin}
V. A. Yudin,
\textit{Multidimensional versions of Paley's inequality,}
(Russian. Russian summary)
Mat. Zametki \textbf{70} (2001),
941--947; translation in
Math. Notes \textbf{70} (2001),
860--865.

\bibitem{Zyg}
A. Zygmund,
\textit{On Fourier coefficients and transforms of functions of two variables},
Studia Math. \textbf{50} (1974), 189--201.

\end{thebibliography}

\end{document}